# THE INSTABILITY OF THE CENTRAL CONFIGURATION CAUSED BY RESONANCE PERTURBATION

Rosaev A.E.

*FGUP NPC Nedra, Svobody 8/38, 150000 Yaroslavl, Russia, rosaev@nedra.ru*

ABSTRACTS: The gravitation interaction between particles is taken into account in central configuration model. After that, the resonance perturbation of the particle of the central configuration by distant satellite is considered in the restricted 3-body problem approach. A shift between simple mean motion resonances and parametric resonance zones is derived. This shift depends upon mass satellite and central configuration properties. Unexpectedly, the resonance's structure depends only on particle's density and number of particles in the central configuration, but not depends on particles mass. Then, the system of few non-interacted planar central configurations is considered. The resonance structure, which follows from parametric instability, is more abundant than for simple mean motion resonance case. The areas of parametric instability restrict stable central configuration, placed in exact commensurability, from both sides. The possible application to the planetary ring structure is discussed.

## 1. INTRODUCTIONS

In some celestial mechanical systems, the motion of a test particle in the gravitational field of a homogeneous ring (or disk) is considered. The perturbation force in this case is calculated as an integral:

$$S_0 = \int_0^{2\pi} \frac{\rho_{surf}(\lambda) d\lambda}{R_0^2 + R^2 - 2RR_0 \cos(\lambda)} \qquad (1.1)$$

where $\rho_{surf}$ - surface density, $R$ and $R_0$ - radial distances.

In the limit, this integral can be replaced by a sum with very large N. At the some cases we can apply central configuration model.

We enter gravitation interaction between particles, modeling central configuration as a system of massive points. Then, the perturbations from distant satellite are included. It is shown, that taking into account interaction between central configuration's particles can sufficiently change resonant structure in system of central configurations.

So, we shall consider central configuration consist of N attractive particles, placed in vertex of regular polygon. Each particle may slightly vary own position and all system rotates with common angular velocity. Set the perturbing body in circular orbit with radius r in central configuration plane. Let the radial perturbation of central configuration to be small, $x \ll R_0$, where $R_0$ - radius of central configuration and $R_0 \ll r$.

We shall consider planar motion of probe particle with mass m in field of attraction other N-1 particles, placed in vertexes of regular polygons and central body with mass M. The equations of motion of the central configuration's particle may be obtained by the calculation of the partial derivatives of perturbation function U (Rosaev, 2003). Evidently, for radial and tangential forces of interaction between particles we have:

$$\ddot{R} - R(\dot{\varphi}+\Omega)^2 = -\frac{GM}{R^2} - \sum_{j=1}^{N-1} Gm_j \frac{2R_0 \sin^2(\alpha_j/2) + x}{\left(x^2 + (2R_0 \sin(\alpha_j/2))^2 (1+x/R_0)\right)^{3/2}}$$

$$R\ddot{\varphi} + 2(\dot{\varphi}+\Omega)\dot{x} = -\sum_{j=1}^{N-1} Gm_j \frac{2R_0 \sin(\alpha_j/2)\cos(\alpha_j/2)}{\left(x^2 + (2R_0 \sin(\alpha_j/2))^2 (1+x/R_0)\right)^{3/2}}$$

(1.2)

where $G$ - gravity constant, $R$ – the central configuration's radius, $\alpha_j$ - angle between particles, and $x$- a distance between the test particle and the central configuration, $R = R_o + x$ – the distance of testing particle from the center, $R_o$ - the distance of $j$-particle from the center, $\Omega$ - angular velocity, $\alpha_j = 2\pi j/N + \varphi$, where φ is possible angular declination from stationary position.

## 2. THEORETICAL MODEL

The equations of planar restricted Hill's problem may be reduced to a Hill equation for normal distance from variation orbit $x$ (Szebehely, 1967):

$$\frac{d^2 x}{dt^2} + \omega^2(t)x = f(t) \tag{2.1}$$

where $f(t)$ – known function of time, $\omega(t)$ - periodic function of time (to be determined).

It means, that area of parametric instability will appear in problem. Note, that variation orbit may be chosen by different way. In particular, it may be two-body problem solution.

At nonperturbed case we have conservation of angular moment $L$. Accounting of perturbation leads to $L$ time dependence. However, $L$ holds constant in average. First of all, the long time effects are interesting, we not interest of central configuration's shape changes. It gives the ability to average satellite perturbation by phase. These expansions confirm it:

$$U_r = -\frac{Gm(R - r\cos(\delta\lambda))}{(R^2 - r^2 - 2Rr\cos(\delta\lambda))^{3/2}} + \left(-\frac{Gm}{(R^2 - r^2 - 2Rr\cos(\delta\lambda))^{3/2}} + \frac{3Gm(R - r\cos(\delta\lambda))^2}{(R^2 - r^2 - 2Rr\cos(\delta\lambda))^{5/2}}\right)x$$

$$U_l = -\frac{GmRr\sin(\delta\lambda)}{(R^2 - r^2 - 2Rr\cos(\delta\lambda))^{3/2}} +$$

$$+ \left(-\frac{Gm}{(R^2 - r^2 - 2Rr\cos(\delta\lambda))^{3/2}} + \frac{3Gm(R - r\cos(\delta\lambda))R}{(R^2 - r^2 - 2Rr\cos(\delta\lambda))^{5/2}}\right) r\sin(\delta\lambda) x$$

(2.2)

For the stability of motion investigation, an averaged equation may be used. An averaged angular momentum is conserved. It means, that variations in mean motion are small, so the problem setting is restricted by conditions:

$$\delta\lambda \equiv \lambda - \lambda_s \approx \delta\varpi t = (\Omega - \omega_s)t \tag{2.3}$$

where $\lambda$ and $\lambda_s$ - mean longitudes, $\Omega$ and $\omega_s$ - mean motions of perturbed and perturbing bodies accordingly, $t$ – time.

Then, we incorporate gravitation interaction between particles, by modeling system of massive points as a central configuration. We shall consider a system to consist of *N* attractive particles, placed close to the vertex of a regular polygon. Each particle may slightly move around stationary position, and the described system rotates with common angular velocity. Set the perturbing body in a circular orbit at radius *r* in the central configuration plane. Let the radial perturbation in system is small, $x \ll R_0$ where $R_0$ - radius of central configuration.

We show, that by taking into account interaction between particles, the behaviour of central configuration at resonance significantly changed.

As it is evident, $\omega$ in (2.1) depend from time:

$$\omega^2 = \omega_0^2 \left(1 + \sum_n h_n \cos n\delta\varpi t\right) \qquad (2.4)$$

where $h_n$ – known coefficient.

Here $\omega_0$ determines the position of centers of unstable zones, and h determines their width $\varepsilon$. Because $\omega_0$ is constant, is possible to restrict largest width (largest $h_n$). Then we can apply this method to each satellite in system. By this way, the construction of the complete structure of instability at central configuration system is possible.

It is sufficient to consider Mathieu equations for the main mean motion resonance cases. At the circular perturbing body orbit:

$$\omega_0^2 = F(R) + \frac{d^2 U_{cfg}}{dR^2} \qquad (2.5)$$

where $F(R)$ – known function of ring central distance $R$, $U_{cfg}$ – central configuration perturbation function. The condition of appearance of parametric resonance and related unstable zones (Landau, 1971):

$$\frac{2\omega_0}{(\Omega - \omega_s)} = n, \qquad n = 1, 2, 3, \ldots \qquad (2.6)$$

where $\Omega$ and $\omega_s$ - mean motions of perturbed and perturbing bodies accordingly. Here *n* – is order of resonance.

As a result, at the system of coaxial central configurations, perturbed by outer satellite, may appear parametric type of instability. The position of resonances is different than in the case where particles gravity is neglected, because the second term in (2.5) appears by central configuration's particles interaction. However, there is one very important effect: the center of instability is shifted relative exact commensurability.

At the case of neglected central configuration particles gravity, this result is in agreement with (Hadjidemetriou, 1982).

The width of resonance is completely determined by satellite perturbation; it rapidly decreases with an resonance order increasing. For the first and second order by k:

$$\varepsilon_1 = \omega_0 \frac{k h_n}{2(k-2)} \qquad \varepsilon_{12} = \omega_0 \frac{k h_n^2}{8(k-2)}$$

$$\varepsilon_1 = \frac{k \omega_s^2 b_n}{2(k-2)\omega_0} \approx \frac{k \omega_s m/M \, b_n}{2(k-2)} \qquad (2.7)$$

where $b_n$ - Laplace coefficients, which easy can be calculated numerically.

## 3. CENTRAL CONFIGURATION AT RESONANCE PERTURBATION

The expansion of central configurations perturbation has a form (Rosaev, 2005):

$$U_r \equiv \frac{dU_{cfg}}{dR} = -\sum_{j}^{n-1}\frac{Gm_j}{R^2(2\sin(\psi_j/2))^2} + \left(\sum_{j=1}^{n-1}\frac{3}{4}\frac{Gm_j}{R^3\sin(\psi_j/2)} - \frac{Gm_j}{R^3((2\sin(\psi_j/2))^3}\right)x +$$

$$\left(\sum_{j=1}^{n-1}\frac{3.75 Gm_j}{R^4((2\sin(\psi_j/2))^3} - \frac{15}{4}\frac{Gm_j}{R^4((8\sin(\psi_j/2))}\right)x^2 + O(x^3) \quad (3.1)$$

where $G$ – gravity constant, $m_j$ – mass of central configuration particle, $\psi_j = \frac{2\pi j}{N} + \phi_0$ - angular distance between the $j^{th}$ – particle and the test particle, $\phi_0$ - small arbitrary angle.

The perturbation function depends on distance relative to the central configuration's axis as shown in Fig.1.

It is possible to rewrite (2.5):

$$\omega^2 = \frac{GM}{R_0^3} + \frac{Gm}{r^3}\sum_{p=2}^{\infty}p(p-1)y^{p-2}P_p(\cos(\delta\lambda)) + \sum_{i=1}^{N-1}\frac{Gm_i}{(2R_0\sin(\pi i/N))^3} \quad (3.2)$$

where N - number of central configuration's particles.

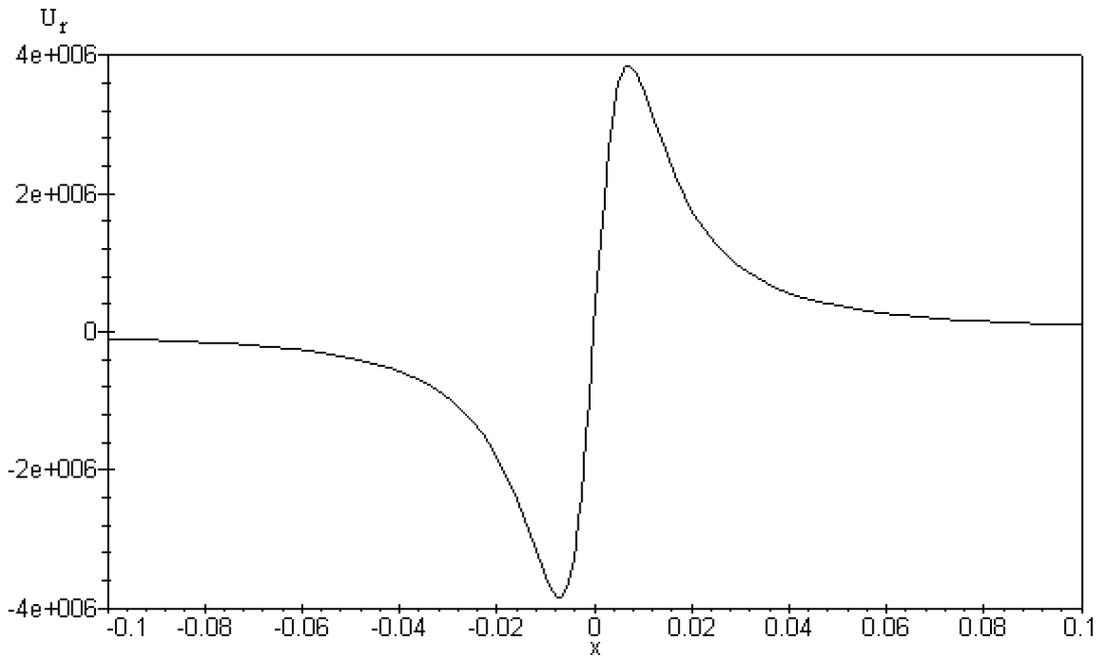

Fig. 1. Central configuration radial perturbation function in dependence of distance from central configuration

## 3.1 SINGLE CENTRAL CONFIGURATIONS

However, there is one very important effect: the center of instability is shift relative exact commensurability. In general case, taking into account approximate expression for mean motions:

$$\Omega^2 \approx \frac{GM}{R^3}, \quad \omega_s^2 \approx \frac{GM}{r^3}, \tag{3.3}$$

resonant condition for ω (2.6) may be rewritten in form:

$$\left(4\alpha + 2\frac{m}{r^3} - \frac{n^2 M}{r^3}\right)\beta^2 + \frac{2n^2 M}{r^{3/2}}\beta + 4M - n^2 M = 0 \quad \beta = R^{3/2} \tag{3.4}$$

Here α - central configuration particles gravity perturbation. It allows one to calculate the positions of centers of instabilities (mass satellite is neglected):

$$R^{3/2} = r^{3/2} \frac{n^2 - 2\sqrt{n^2 + (n^2 - 4)\alpha r^3 / M}}{n^2 - 4\alpha r^3 / M} \tag{3.5}$$

There are two cases at strong interaction between particles. The sign of α is determined by expansion (3.1). In the case of single central configuration α is positive, and we have instability relative small oscillations – a very narrow gap, shifted toward the planet.

In case of α positive, it is easy to calculate:

$$\sum_{i=1}^{N-1} \frac{Gm_i}{(2R_0 \sin(\pi i / N))^3} \approx 0.300512625 \frac{Gm_i N^3}{\pi^3 R_0^3} = 1.258784 G\rho\sigma^3 = G\alpha \tag{3.6}$$

where $\sigma = r_k N/(\pi R_0)$ - where σ, ρ, and $r_k$ - central configuration's «coefficient of filling», density and size of central configuration's particles accordingly. The case σ=1, for example, corresponds to the model of maximal density, where particles in central configuration touch one another.

At large amplitude of perturbation, the parametric instability may be (more) significant, and according unstable zones are shifted from exact commensurability. It is very important, that this condition not depend on the central configuration (or central configuration's particle) mass, only from degree of filling the central configuration by particles. *It means, that even for low-massive central configurations, mutual gravitation interaction between particles can play a significant role in central configuration particles dynamics, at least, at small area close to resonance and close to central configuration.*

The described shift depends on particle size and order of resonance (Fig 2 -3).

For stable motion (provide stability of motion) central configuration must rotate with velocity, higher, than keplerian (Rosaev, 2003 ):

$$\Omega = 1/\sqrt{R_0} \left( \frac{GM}{R_0^2} + \sum_{j=1}^{N-1} Gm \frac{1}{4R_0^2 \sin(\alpha_j/2)} \right)^{1/2} = \Omega_0 \left( 1 + \frac{m}{4M} \sum_{j=1}^{N-1} \sin^{-1}(\alpha_j/2) \right)^{1/2} \tag{3.7}$$

This fact effects shift the central configuration away from central mass. However, in most cases it is negligible.

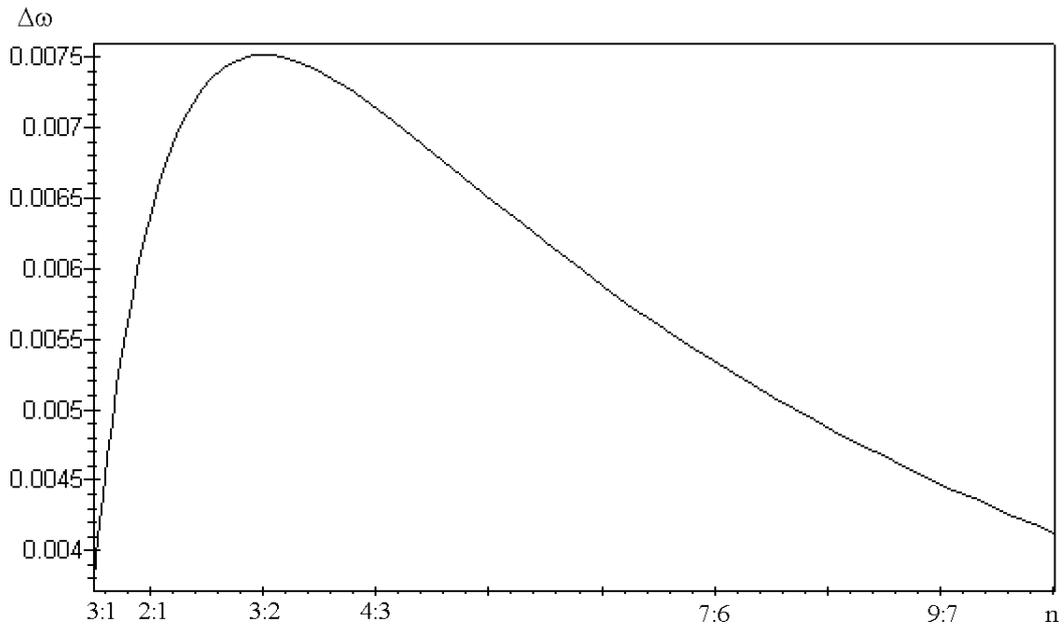

Fig. 2. Shift in instable zone central distance due to mutual attraction of central configurations particles ($\alpha=0.1$, $m=0$) in dependence of order of resonance

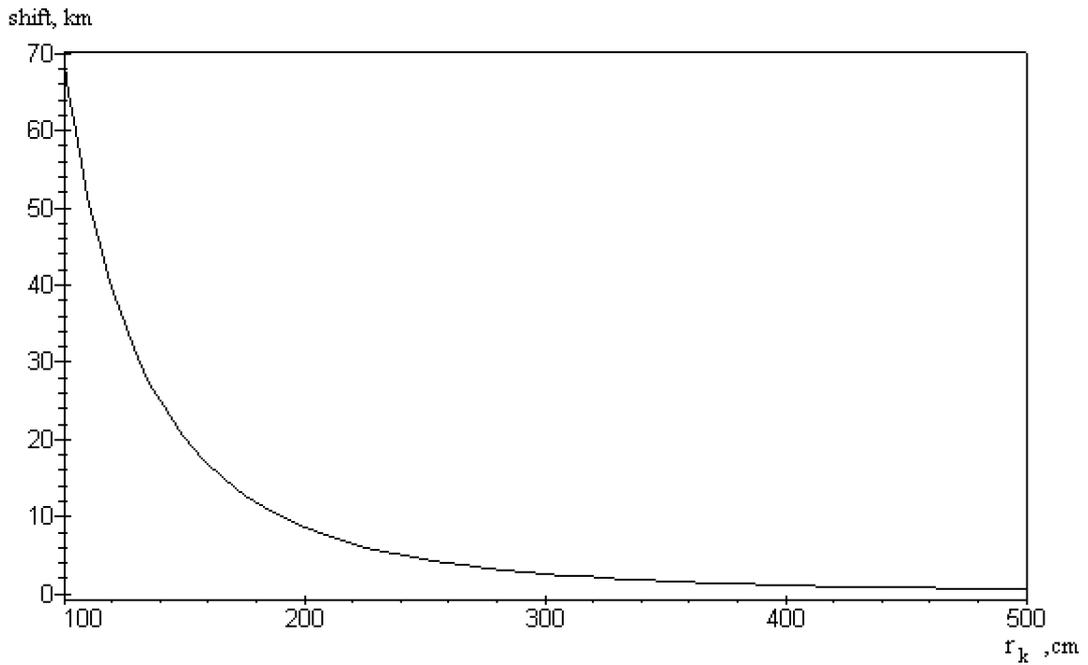

Fig 3. The shift of resonances due to particles gravity in dependence of particle size

## 3.2. A SYSTEM OF CENTRAL CONFIGURATIONS

Let us consider the system with 2N masses m, which form two regular polygons, the vertex of each lying on the radius the circle described (Fig.4). Our target is to find areas of instability for the second central configuration, perturbed by the satellite on resonant orbit and by the first central configuration.

It may be assumed, that in case of insignificant distance $d$ between the rings, the such system may exist in reality, rotating with average angular velocity:

$$\Omega_k = (\Omega_o + \Omega_i)/2 \tag{3.8}$$

The requirement of common angular velocity at a small $d$ leads to the following conditions:

$$\frac{3M}{m_k} d_k / R_0 \approx \sum_{j=1}^{N-1} \frac{1}{4|\sin(\alpha_j/2)|} + \frac{R_0^3}{d_k^3} \tag{3.9}$$

This is a remarkable result: $m_k$ increase with $k$. It means, that the mass of the ring's particles should increase outwards to preserve the common angular velocity. Initial $d$ correlates with libration point distances, finite $d_n$ are determined by $M/m$ ratio. The angular moment is constant in described system (L= $R^2 \Omega_0$), and after linearization, equation (10) can be re-written in the form:

$$\ddot{x}_k + (3L^2/R^4 - 2GM/R^3)x_k = B + \frac{d\widetilde{F}_k}{dx} x_k - \frac{d\widetilde{F}_s}{dx} x_k - 2\frac{Gm}{d^3} x_k + \ldots$$

$$R\ddot{\varphi}_k + \frac{df_k}{d\varphi_k}\varphi_k - 2\frac{Gm}{d^3}\psi = 0; \qquad \psi/\varphi \approx R/d \tag{3.10}$$

B= $B_k$ -$B_s$

The expressions (3.9) and (3.10) give the following conditions for radial and tangential stability ($d$ is distance between the rings):

$$\Omega_{1o}^2 - \frac{2Gm_1}{d^3} + Gs(m_o - m_i + 3d/Rm_o) > 0$$

$$\frac{d(f_i + f_o)}{d\varphi} - \frac{2Gm_1}{d^3} R < 0$$

$$\tag{3.11}$$

$$\Omega_{1i}^2 - \frac{2Gm_o}{d^3} + Gs(m_i - m_o + 3d/Rm_o) > 0$$

$$\frac{d(f_{i1} + f_o)}{d\varphi} - \frac{2Gm_o}{d^3} R < 0$$

where:

$$\Omega_{1k}^2 = \Omega_k^2 + \frac{3G(m_s - m_k)}{8R^2} \sum_{j=1}^{N-1} \frac{1}{|\sin(\pi j/N)|} \quad (3.12)$$

Here, it is allowed for Ro = Ri + d = R + d. It is evident, that only radial perturbations are important for the stability at small $d$. As a whole, the central mass $\Omega_1 \sim \Omega_0$. In this case, taking into account, that $d>0$, except for $\Omega_0$ from (3.11), the conditions of linear stability lead to:

$$1 - 3\,d/R < m_i/m_o < 1 \quad (3.13)$$

It means, that stability is possible only when the outer mass exceeds the inner one. At different limit M→0:

$$\Omega_{1k}^2 = \frac{3G(m_s - m_k)}{8R^2} \sum_{j=1}^{N-1} \frac{1}{|\sin(\pi j/N)|} \quad (3.14)$$

Evidently, at each fixed $mi$ and $mo$, a limit of stability exist depending on $d$. At small $d$ collisions between $mi$ and $mo$ is inevitable, and two rings join to one. The limit value of $d$ can be obtained from the expression:

$$d > R/3\left[(m_o - m_i)\Omega_k^2 + s(m_o^2 - m_i^2)\right]/(m_o s(m_o + m_i)) \quad (3.15)$$

with

$$s = A_k/Gm_k > 0 \quad (3.16)$$

Finally, the conditions of stable configuration existence, rotating with the same angular velocity, can be obtained by substitution of (17) and (21) to (20). The dimension of such configuration and the mass distribution within the described structure are given in the equations:

$$d > \frac{R(m_o A_k/Gm_k(m_o + m_i))/3}{\left[(m_o - m_i)\dfrac{3G(m_s - m_k)}{8R^2}\sum_{j=1}^{N-1}\dfrac{1}{|\sin(\pi j/N)|} + A_k/Gm_k(m_o^2 - m_i^2)\right]}$$

$$(3.17)$$

$$A_k = -\frac{Gm_k}{R^3}\sum_{j=1}^{N-1}\left[\frac{1}{|2\sin(\alpha_j/2)|^3} - \frac{3}{|8\sin(\alpha_j/2)|}\right]$$

In (Seidov Z.F., 1990) it is assumed, that the system of two central configurations can be stable in case shown in Fig1a. According to our considerations, we can expect, that the described configuration (Fig. 4, right) is stable in linear approximation. In contrast, the configuration in Fig.4 (left) is probably loosing stability due to interaction of two rings.

It is possible to apply piece-wise linearization of the central configuration potential (Fig. 5).

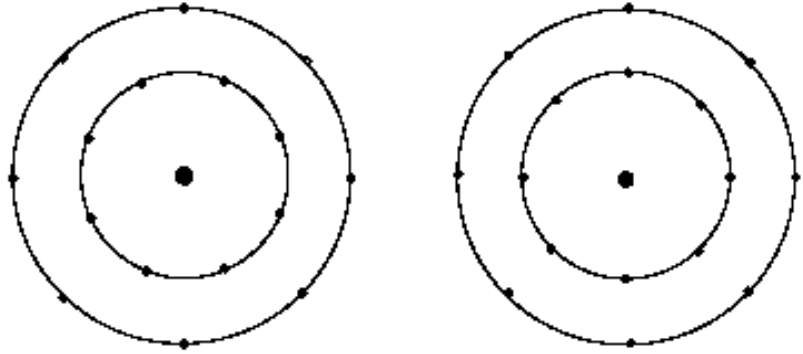

Fig. 4. Non-collinear and collinear 2N + 1 body configurations.

In case of α negative, we have oscillation with large amplitude: it is sufficient, if the particle will be out first linear zone to remove it far from it's initial position. So, we have a broad area of instability, shifted away from central planet relative exact commensurability. In result, we have two unstable gaps at each resonance.

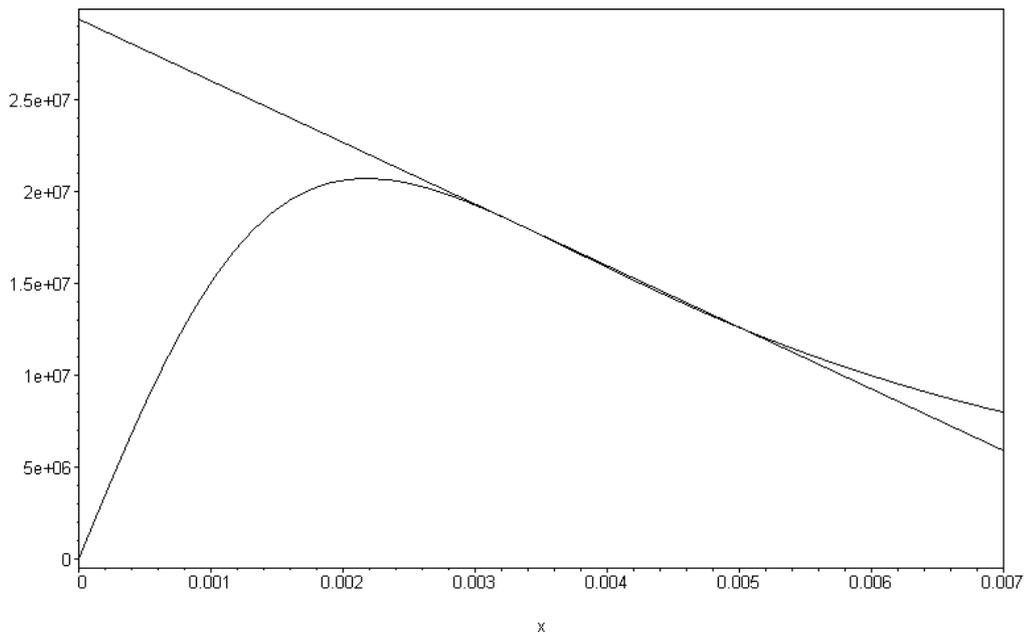

Fig.5. Linearization at negative α

For this case, the position of unstable zones is calculated by:

$$R = r \left[ \frac{n^2 - 2\sqrt{n^2 + (n^2 - 4)\alpha_1 r^3 / M}}{n^2 - 4\alpha_1 r^3 / M} \right]^{2/3} \pm \Delta x \qquad (3.18)$$

where

$$\alpha_1 \approx -0.5\alpha, \qquad \Delta x \approx \frac{2\pi R}{N} \text{ – determined by area of linear motion}$$

At the most simplest case σ=const we have two time more resonances in central configuration region, than in followed from 3-body problem (Fig. 6). This situation become possible, when distance from central mass overcome critical value $R_c$ and coefficient of fullness σ sufficiently large. Note, that at fixed surface density, the minimal possible size of particles, take contribution in gravity, exist and determined by σ.

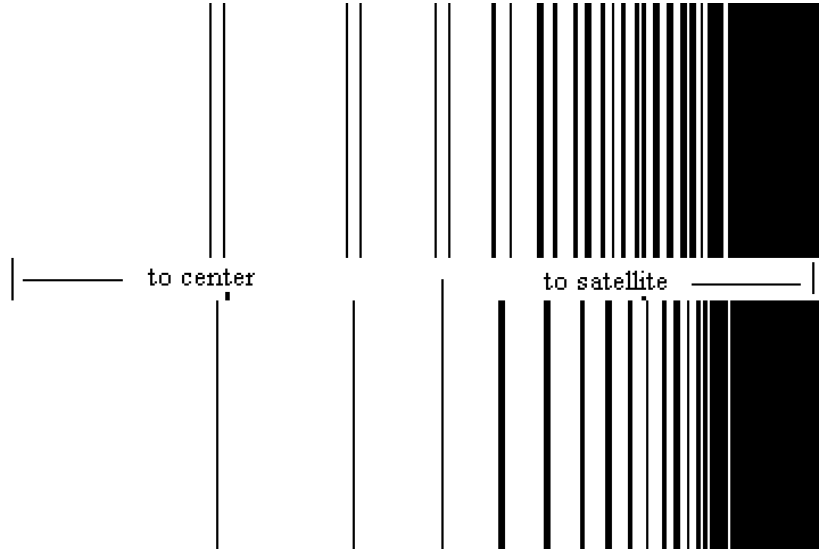

Fig.6. Resonance structure at α=0.5

## 4.APPLICATIONS

As a possible application of model, the system of planetary rings maybe considered. Recently, the model of the central configurations successfully applied to the real planetary rings by Meyer K.R. and Schmidt D.S. (Meyer, Schmidt, 1993). By using for surface density $\rho_{surf}=60$ g/cm² (Stone *et al.,* 1982), $r_k = 1$ m, $\rho = 1$ g/cm³, the estimations of number of particles N, σ and α gives:

$$N = \frac{3R\rho_{surf}}{2\rho r_k^2} \approx 9*10^7 \quad \sigma = \frac{\rho_{surf}}{4\rho r_k} = 0.15 \div 0.3$$
$$\alpha = 1.25\rho\sigma^3 = 1.25\rho\left(\frac{Nr_k}{2\pi R}\right)^3 \approx 4 \div 12$$
(4.1)

The parametric instability effect most significant at strong interaction between particles ( $M/R^3 < 1.25\rho\sigma^3$ ) The according limits for real planetary ring system are given in table1. As a result, possible effect on main resonances at saturnian ring is shown in table 2.

Table 1

Possible conditions of parametric instability in planetary rings

| Ring system | Range of distance (thous. km) | Min distance for instability (thous. km) | Coefficient of filling, σ |
|---|---|---|---|
| Jupiter | 127-130 | 103 | 0.80 |
| Saturn | 70-140 | 69 | 0.49-0.98 |
| Uran | 41-51 | 36.5 | 0.72-0.90 |
| Neptun | 41-62 | 38.2 | 0.63-0.95 |

Table 2

The positions of the main parametric instability areas within Saturn B ring.

| N | Resonance | Distance, $*10^5$ km | | | Satellite |
|---|---|---|---|---|---|
| | | ($\rho_{surf}=0$ g/cm$^2$, $\rho = 0$ g/cm$^3$) | ($\rho_{surf}=60$ g/cm$^2$, $r_k = 0.5$ m, $\rho = 1$ g/cm$^3$) | | |
| 1 | 5:3 | 0.95040189 | 0.94784459 | 0.95301368 | Pan |
| 2 | 2:1 | 0.95421381 | 0.95153764 | 0.95694358 | Epimeteus |
| 3 | 5:3 | 0.97914159 | 0.97650695 | 0.98183236 | Atlas |
| 4 | 5:3 | 0.99130616 | 0.98863880 | 0.99403036 | Prometeus |
| 5 | 5:3 | 1.00802356 | 1.00531121 | 1.01079370 | Pandora |
| 6 | 3:2 | 1.08137339 | 1.07867963 | 1.08412544 | Pandora |
| 7 | 3:2 | 1.01955882 | 1.01701904 | 1.02215356 | Pan |
| 8 | 3:2 | 1.05038979 | 1.04777321 | 1.05306299 | Atlas |
| 9 | 3:2 | 1.06343953 | 1.06079044 | 1.06614595 | Prometeus |
| 10 | 5:3 | 1.07753949 | 1.07464008 | 1.08050066 | Epimeteus |
| 11 | 3:1 | 1.14428081 | 1.14142243 | 1.14718455 | Encelad |
| 12 | 3:2 | 1.15594770 | 1.15306817 | 1.15888955 | Epimeteus |
| 13 | 2:1 | 1.16870277 | 1.16542505 | 1.17204614 | Mimas |

## 5. CONCLUSIONS

The resonance perturbation of the system of the planar central configurations by the distant satellite is considered. Our method consists of linearization equation of motion and reduction of them to Hill or Mathieu equation to study parametric instability in system. All results are valid as for system of non-interacted central configurations as well as for single central configuration.

The resonance structure, followed from parametric instability, at least two times more abundant than for simple mean motion resonance case.

A shift between simple mean motion resonances and parametric resonance zones is derived. This shift depends on mass satellite and central configuration properties. Unexpectedly, the resonance's structure depends only on the particle's density, but does not depend on the particle's (or central configuration's) mass. It is shown, that mutual interaction between particles of ring causes a significant effect on resonance structure.

There are 3 main results of mutual gravity between particles taking into account: 1) shift instable gaps out of center from exact commensurability, 2) more abundant resonance structure, 3) possibility single-side shepherding, when one satellite can bound central configuration from both sides.

The results may be applied to Saturn ring system. It is shown, that by varying the ring's particles properties, is possible to explain observed shift between actual ring distances and exact commensurability.